\documentclass{elsart3-1}

\usepackage{latexsym}
\usepackage{amsmath}
\usepackage{amsfonts}
\usepackage{amssymb}
\usepackage{graphpap}
\usepackage{epsfig}
\usepackage{amscd}
\usepackage{enumerate}
\usepackage{eucal} 

\newcommand{\RR}{{\mathbb{R}}}
\newcommand{\NN}{{\mathbb{N}}}
\newcommand{\ZZ}{{\mathbb{Z}}}
\newcommand{\CC}{{\mathbb{C}}}

\newcommand{\eps}{\varepsilon}

\newcommand{\ep}{\hfill$\Box$ \vskip 0.08in}

\usepackage[english,francais]{babel}

\newtheorem{theorem}{Theorem}[section]

\newtheorem{e-proposition}[theorem]{Proposition}

\newtheorem{e-definition}[theorem]{Definition\rm}

\newtheorem{theoreme}{Th\'eor\`eme}[section]

\newtheorem{proposition}[theoreme]{Proposition}

\renewcommand{\theequation}{\arabic{equation}}
\def\og{\leavevmode\raise.3ex\hbox{$\scriptscriptstyle\langle\!\langle$~}}
\def\fg{\leavevmode\raise.3ex\hbox{~$\!\scriptscriptstyle\,\rangle\!\rangle$}}

\journal{}
\begin{document}
\centerline{}
\begin{frontmatter}


\selectlanguage{english}
\title{Boundedness of the square function and rectifiability}


\selectlanguage{english}
\author[SM]{Svitlana Mayboroda},
\ead{svitlana@math.purdue.edu}
\author[AV]{Alexander Volberg}
\ead{volberg@math.msu.edu, A.Volberg@ed.ac.uk}

\address[SM]{Department of Mathematics, Purdue University, 150 N. University Street, West Lafayette, IN 47907-2067, USA}
\address[AV]{Department of Mathematics, Michigan State University, East Lansing, MI 48824, USA}


\medskip
\begin{center}
\end{center}

\begin{abstract}
\selectlanguage{english}
Following a recent paper \cite{Tolsa1} we show that the finiteness of square function associated with the Riesz transforms with respect to Hausdorff measure $H^n$ ($n$ is interger) on a set $E$ implies that $E$ is rectifiable.


\vskip 0.5\baselineskip

\selectlanguage{francais}

\end{abstract}
\end{frontmatter}



\selectlanguage{english}
\section{Introduction}
\label{}


\font\tretamajka=cmbsy10 at 12pt

For a Borel measure $\mu$ in $\RR^d$ and $s\in (0,d]$ the $s$-Riesz transform of $\mu$ is defined as
\begin{equation}\label{eq1.1}
R^s \mu(x):=\int\frac{x-y}{|x-y|^{s+1}}\,d\mu(y), \qquad x\not\in {\rm supp}\,\mu,
\end{equation}

\noindent and the truncated Riesz transform is given by
\begin{equation}\label{eq1.2}
R^s_{\eps}\, \mu(x):=\int_{|x-y|>\eps}\frac{x-y}{|x-y|^{s+1}}\,d\mu(y), \quad R^s_{\eps,\eta}\, \mu(x):=R^s_{\eta}\, \mu(x)-R^s_{\eps}\, \mu(x),
\end{equation}

\noindent where $x\in\RR^d,$ $\eta>\eps>0$.

One says that the set $E\subset \RR^d$ is $n$-rectifiable, $n\in\NN$, if there is a countable family of $n$-dimensional $C^1$ submanifolds  $\{M_i\}_{i\geq 1}$ such that $H^n\left(E\setminus \cup_i M_i\right)=0$, where $H^n$ stands for the $n$-dimensional Hausdorff measure.

The main result of this note is as follows. 

\begin{theorem}\label{tSqF} Let $\mu$ be a finite Radon measure in $\RR^d$. Assume that for some $s\in(0,d]$ there is a set $E\subset \RR^d$ with the property that
\begin{equation}\label{eq1.2.3}
0<\theta^{s,\ast}_{\mu}(x)<\infty\quad\mbox{for all}\quad x\in E,
\end{equation}

\noindent and the square function
\begin{equation}\label{eq1.3}
S\mu(x):=\Bigg(\int_0^\infty \left|R^{s}_{t,2t} \mu(x)\right|^2\,\frac{dt}{t}\Bigg)^{1/2}<\infty \quad\mbox{for all}\quad x\in E.
\end{equation}

\noindent  Then $s$ is an integer and $E$ is $s$-rectifiable.

In particular, if $E$ is a compact set in $\RR^d$ with $0<H^s(E)<\infty$ and \eqref{eq1.3} is satisfied, then $s$ is an integer and $E$ is $s$-rectifiable.
 \end{theorem}

Here and below, the upper $s$-dimensional density of $\mu$ at $x$ is given by  
$\theta^{s,\ast}_{\mu}(x):=\limsup_{r\to 0}\frac{\mu(B(x,r))}{r^s}$, where $B(x,r)$ is the ball of radius $r>0$ centered at $x\in \RR^m$. 

In  \cite{DS1}, \cite{DS2} David and Semmes showed, under certain assumptions on the measure $\mu$, that the $L^2$ boundedness of {\it all}  Calder\'on-Zygmund singular integral operators implies that $s$ is an integer and $\mu$ is $s$-rectifiable. Our main result assures that it is, in fact, sufficient to assume  pointwise boundedness of a single operator, namely, the square function, in order to arrive to the same conclusion.  Alternatively, one could view  \eqref{eq1.3} {\it almost} as a condition
\begin{equation}\label{eq1.4}
{\mathbb{E}}\,\Bigl|\sum_{k\in\ZZ} \eps_k R^s_{2^{-k},2^{-k+1}}\, \mu(x)\Bigr|<\infty \qquad \mu - {\rm a.e.}\,x\in\RR^m,
\end{equation}

\noindent where $\eps_k$ are independent random variables taking the values $-1$ and $1$ with probability $1/2$ each. Therefore, roughly speaking,  the boundedness of the singular integrals of the type $\sum_{k=0}^{\infty} \eps_k R^s_{2^{-k},2^{-k+1}}\, \mu(x)$ already guarantees that $\mu$ is $s$-rectifiable, and the corresponding  class of operators is much smaller than that of David and Semmes.

Let us point out that we have already proved the fact that \eqref{eq1.2.3}--\eqref{eq1.3} imply that $s$ is an integer in \cite{MV1}. The present paper, concentrating on the issue of rectifiability, is a sequel to the aforementioned work. 

Ultimately, one would like to show that the conclusion of Theorem~\ref{tSqF} holds purely under the assumption that the Riesz transforms $R_\eps^s\mu$ are bounded in $L^2$ uniformly in $\eps>0$, or under the condition 
\begin{equation}\label{eq1.2.2}
\sup_{\eps>0} |R_{\eps}^s\,\mu(x)|<\infty\qquad \mu - {\rm a.e.}\,x\in\RR^d.
\end{equation}
\noindent This is a long-standing open problem, known as a conjecture of Guy David. At the moment, it has been resolved only for $0<s\leq 1$ (\cite{To2},  \cite{Pr2}). The proof heavily relies  on the curvature of measure estimates from \cite{MeV}, which are decisively restricted to lower dimensions. However, in the presence of Theorem~\ref{tSqF}, the problem of Guy David in all dimensions is essentially reduced  to the passage from \eqref{eq1.2.2} to \eqref{eq1.3}.

Our methodology builds on recent advances in \cite{Tolsa1}, where the author proved that the existence of $p.v. R^n\mu(x):=\lim_{\eps \to 0} R^n_\eps \mu(x)$ for $\mu$-a.e. $x\in\RR^d$ entails that $\mu$ is $n$-rectifiable. Analogously to the argument in \cite{Tolsa1}, one of the main ingredients in our proof is an estimate for the $L^2$ norm of the square function, on the graph of a Lipschitz function $A$, in terms of $\|\nabla A\|_2$. It compensates for the lack of the curvature of measure techniques. Having this at hand, the matters are further reduced to the L\'eger's construction of Lipschitz graphs \cite{Le}. In the present note  we outline the main stages of the proof, and the details will appear elsewhere.

\section{$L^2$ estimates on a Lipschitz graph}

Let us start by introducing some extra notation. As we already mentioned, the case of the non-integer dimension was treated in \cite{MV1}, and here we focus on $s=n$, $n\in\NN$, $n\leq d$.
Let $\Pi$ stand for the projection $(x_1,...,x_n,...,x_d)\to (x_1,...,x_n,0,...,0)$, $\Pi^{\perp}:=I-\Pi$.

Furthermore, let $R^{j} \mu(x)$, $k\in\ZZ$, $x\in\RR^d$, be a version of doubly truncated Riesz transform, defined, for example, as follows. If $\psi^0$ is a non-increasing radial $C^\infty$ function with $\chi{(B(0,1/2))}\leq \psi^0\leq \chi(B(0,2))$ we set $\psi_j(z):=\psi^{0}(2^jz)-\psi^0(2^{j+1}z)$, $j\in\ZZ$, $z\in\RR^d$, so that each $\psi_j$ is non-negative,  ${\rm supp}\,\psi_j\subset B(0,2^{-j+1})\setminus B(0,2^{-j-2})$ and $\sum_{j\in\ZZ}\psi_j=1$ in $\RR^d\setminus\{0\}$.  Then $R^{j} \mu$ is an operator defined analogously to \eqref{eq1.1} with the kernel given by $\psi_j(x-y)\,\frac{x-y}{|x-y|^{s+1}}$. 


It is not hard to show that the condition \eqref{eq1.3} implies that the discrete version of the square function
\begin{equation}\label{eq1.3d}
\widetilde S\mu(x):=\Bigg(\sum_{j\in\ZZ} \left|R^{j} \mu(x)\right|^2\Bigg)^{1/2}<\infty \quad\mbox{for all}\quad x\in E.
\end{equation}

In the course of the proof we will also employ the notation $\widetilde S^\perp \mu:=(\sum_{j\in\ZZ} \left|\Pi^\perp(R^{j} \mu)\right|^2)^{1/2}$ and $S^\perp\mu=(\int_0^\infty \left|\Pi^\perp(R^{n}_{t,2t} \mu)\right|^2\,\frac{dt}{t})^{1/2}$.

Our reasoning follows closely to \cite{Tolsa1}, replacing $R$'s by $S$'s, and one of the cornerstones of  our argument is the following $L^2$ bound for the square function on a Lipschitz graph.

\begin{theorem}\label{tL2} Let $\Gamma$ be an $n$-dimensional Lipschitz graph $\Gamma=\{(x,y)\in \RR^n\times \RR^{d-n}:\,y=A(x)\}$ and let $d\mu=g(z)dH^n_{|\Gamma}(z)$ with $C_1^{-1}\leq g(z)\leq C_1$ for some $C_1>0$ and all $z\in\Gamma$. Assume further that $A$ is compactly supported, $\|g-1\|_2\leq C_2 \|\nabla A\|_2$, and $\|\nabla A\|_{\infty}\leq \eps_0$. Then for $\eps_0=\eps_0(C_2)$ sufficiently small 
\begin{equation} \label{eqL2}
\|S^{\perp}\mu\|_{L^2(\mu)}\approx \|S\mu\|_{L^2(\mu)}\approx \|\nabla A\|_2.
\end{equation}
\end{theorem}

\noindent {\it Sketch of proof}.\, The upper estimate $\|S\mu\|_{L^2(\mu)}\leq C\|\nabla A\|_2$ does not require the smallness of $\eps_0$,
and follows directly by combining $(5.10)$ in \cite{TolsaUnif} with the argument of Lemma~3.1 in \cite{Tolsa1}. The lower bound, $\|S^\perp\mu\|_{L^2(\mu)}\geq C\|\nabla A\|_2$, is more involved.

Observe that $\Pi_{\sharp}\mu=\rho(x)\,dx$ with $\rho(x)=g(x)J\widetilde A(x),$ where $\widetilde A(x)=(x,A(x))$, $x\in\RR^n$, and $J\widetilde A$ stands for the $n$-dimensional Jacobian of $\widetilde A$. Take now $\mu_0$ with ${\rm supp}\,\mu_0\subset\Gamma$ such that $\Pi_{\sharp} \mu_0=dx$. Then $d\mu(x)-d\mu_0(x)=(h(x)-1)\,d\mu_0 (x)$, $h=\rho(\Pi(x))$, and one can show that  $\|\rho-1\|_2\leq C\|\nabla A\|_2$, so that $\|h-1\|_2\leq C\|\nabla A\|_2$. 

Carefully tracking the argument in Chapter~5 of \cite{Tolsa1} we deduce that for such a measure $\mu_0$ 
\begin{equation}\label{eqL2.1}
\|S^\perp\mu_0\|_{L^2(\mu_0)}^2\approx \sum_{j\in\ZZ}\|\Pi^\perp(R^j\mu_0)\|_{L^2(\mu_0)}^2\approx \|\nabla A\|_2^2,
\end{equation}

\noindent hence, it remains to estimate the difference between $\|S^\perp\mu_0\|_{L^2(\mu_0)}$ and $\|S^\perp\mu\|_{L^2(\mu)}$. However,
\begin{equation}\label{eqL2.2}
\left|\|S^\perp\mu_0\|_{L^2(\mu_0)}-\|S^\perp\mu\|_{L^2(\mu_0)}\right|\leq \left\|S^\perp(\mu-\mu_0)\right\|_{L^2(\mu_0)}=\left\|S^\perp((h-1)\,d\mu_0)\right\|_{L^2(\mu_0)}.
\end{equation}

 Furthermore, when $\Pi_{\sharp} \mu_0=dx$, one can directly show that the operator $S^{\perp}$ is bounded in $L^2(\mu_0)$, with the norm controlled by $C\|\nabla A\|_{\infty}$.  Hence, the right-hand side of \eqref{eqL2.2} is bounded from above by $
C\|\nabla A\|_{\infty}\|h-1\|_{L^2(\mu_0)}\leq \eps_0\|\nabla A\|_2.$
\noindent Then, if $\eps_0$ is sufficiently small, \eqref{eqL2.1}--\eqref{eqL2.2} lead to $\left\|S^\perp\mu\right\|_{L^2(\mu_0)}\approx \|\nabla A\|_2$, which implies that $\left\|S^\perp\mu\right\|_{L^2(\mu)}\approx \|\nabla A\|_2$ since $g(x)\approx h(x)\approx 1$.\ep

\section{From $E$ to the construction of Lipschitz graphs}

Theorem~\ref{tSqF} can be reduced to the following Proposition mimicking the main Lemma of \cite{Tolsa1}.
\begin{proposition}\label{MainLemma} Let $\mu$ be a finite Borel measure in $\RR^d$. Assume that there is a closed ball $B_0=\overline{B} (x_0,r_0)$ and a compact set $F\subset 10B_0$ with $x_0\in F$ such that for some positive constants $M_1,M_2,\delta_1,\delta_2$
\begin{enumerate}
\item $\mu(8B_0)=c_n8^nr_0^n$ and $\mu(10 B_0\setminus F)\leq\delta_1\mu(B_0)$,
\item $\mu(B(x,r))\leq M_1r^n $ for all $x\in F$, $r>0$, and $\mu(B(x,r))\leq c_n(1+\delta_1)r^n$ for all $x\in F$ and $0<r\leq 100 r_0$,
\item $\|S\|_{L^2(\mu |_F),L^2(\mu |_F)}\leq M_2$,
\item $|\widetilde R_{\eps,2\eps} \mu(x)|+|\widehat R_{\eps,2\eps} \mu(x)|\leq \delta_2$ for all $x\in F$ and $0<\eps<\delta_2^{-2}r_0$.
\end{enumerate}

If $\delta_1=\delta_1(M_2)$ and $\delta_2=\delta_2(M_1,M_2)$ are small enough, then there exists an $n$-dimensional Lipschitz graph $\Gamma$ with the property that 
\begin{equation}\label{eqMainLemma1}
\mu(\Gamma\cap F\cap B_0)\geq {\textstyle{\frac{9}{10}}}\,c_nr_0^n.
\end{equation}
\end{proposition}

Here, 
\begin{equation}\label{eq1.2tilde}
\widetilde R_{\eps}\, \mu(x):=\int\frac{x-y}{(\eps^2+|x-y|)^{\frac{n+1}{2}}}\,d\mu(y), \quad \widehat R_{\eps}\mu(x)=\int \psi\left(\frac{x-y}{\eps}\right)\frac{x-y}{|x-y|^{n+1}}\,d\mu(y),
\end{equation}

\noindent for some $C^\infty$ function $\psi$ with $\chi_{\RR^d\setminus B(0,1)}\leq \psi\leq \chi_{\RR^d\setminus B(0,1/2)}$, and $\widetilde R_{\eps,2\eps}\, \mu=\widetilde R_{2\eps}\, \mu-\widetilde R_{\eps}\, \mu,$ $\widehat R_{\eps,2\eps}\, \mu=\widehat R_{2\eps}\, \mu-\widehat R_{\eps}\, \mu$ , $\eps>0$.

The reduction from Theorem~\ref{tSqF} to Proposition~\ref{MainLemma} follows the general lines of the argument in Chapter~7 of \cite{Tolsa1}. An important new element  is the need  to extract the sets satisfying the condition $(iii)$ above. This is also interesting on its own right,  and is stated in the following Proposition.

\begin{proposition}\label{From_pw_to_L2} Let $\mu$ be a positive finite Radon measure on $\RR^d$, satisfying $0<\theta_\mu^{n,\ast}(x)<\infty$ for $\mu$-a.e. $x\in\RR^d$, $n\in\NN$, $n\leq d$, and such that $S\mu(x)<\infty$ for $\mu$-a.e. $x\in\RR^d$. Then for every $\delta>0$ there is a compact set $E\subset \RR^d$ such that $S$ is bounded in $L^2(\mu |_E)$ and $\mu(\RR^d\setminus E)\leq \delta$.
\end{proposition}

\noindent {\it Sketch of Proof}.\,
Following the techniques in \cite{VBook} and \cite{To2}, we prove Proposition~\ref{From_pw_to_L2} via reduction to the operators with suppressed kernels. To this end, let us denote by $S_\eps$ the square function defined analogously to \eqref{eq1.3} with the integration over $t\in(\eps,\infty)$. Also, let $S^\theta$ be an analogue of \eqref{eq1.3} with the integrand given by $R^\theta_{t,2t}$, where $\theta$ is a non-negative Lipschitz function with Lipschitz constant less than or equal to 1, and 
\begin{equation}\label{eqDefSup}
R^\theta_{t,2t}\mu(x)=\int_{t\leq |x-y|\leq 2t}k_\theta(x,y)\,d\mu(y), \quad k_\theta(x,y)=k(x,y)\,\frac{1}{1+k^2(x,y)\theta(x)\theta(y)},
\end{equation}

\noindent with $k(x,y)=\frac{x-y}{|x-y|^{n+1}}$ being the kernel of the Riesz transform. 

First, we show that the square function with the suppressed kernel $S^\theta \mu$  is controlled by $S_\eps \mu$, in the following sense. If for some $x\in\RR^d$ and  $r_0\geq 0$
\begin{equation}\label{eqLemmaSup}
\mu(B(x,r))\leq C_0r^n\quad\forall\,\,r\geq r_0,\quad\mbox{and}\quad |S_\eps\mu(x)|\leq C_1\quad \forall\,\,\eps\geq r_0, \end{equation}

\noindent then $|S_\theta \mu(x)|\leq C_2(C_0,C_1)$, provided that $\theta(x)\geq r_0$. Next, for every $\delta>0$ one can find an open set $H$ containing all non-Ahlfors balls (that is, all balls with the property $\mu(B)>Mr^n$ for some constant $M=M(\delta)$) with $\mu(H)\leq \delta$  and such that $|S_\eps\mu(x)|\leq C(\delta)$ for every $\eps>{\rm dist}\,\{x,\RR^d\setminus H\}$, for all $x\in\RR^d$. According to the above, this, in particular, implies that $|S_\theta \mu(x)|\leq C'(\delta)$ for every Lipschitz function $\theta(x)\geq {\rm dist}\,\{x,\RR^d\setminus H\}$ with the Lipschitz constant 1, for all $x\in\RR^d$. We claim that for every such $\theta$ 
\begin{equation}\label{eqTb}
|S_\theta\mu(x)|\leq C \,\,\mbox{for every}\,\,x\in \RR^d\quad\Longrightarrow \quad S_\theta:L^2(\mu)\to L^2(\mu).
\end{equation}

Indeed, let $\overline S_\theta$ be the mapping from $f:\RR^d\to L^2(dt)$ to 
\begin{equation}\label{eqTb1}
\overline S_\theta f(x,t):=\frac{1}{\sqrt t}\int_{t\leq |x-y|\leq 2t} k_\theta(x,y)f(y,t)d\mu(y).
\end{equation}

\noindent Then one can employ the non-homogeneous $T1$ theorem for vector-valued kernels ($Tb$ Theorem~4 in \cite{Hyt}), based, in turn, on the earlier developments in \cite{NTV-Acta}. We omit in this note the detailed verification of the conditions of the Theorem. The most important part is that
the bound $S_\theta\mu(x)\leq C$, $x\in\RR^d$, amounts to $\overline S_\theta ({\bf1}) \in L^\infty (\mu, L^2(dt))$, where ${\bf 1}$ is the identity on $\RR^d\times (0,\infty)$. This, in a sense, plays a role of the $T1$ condition, and allows to conclude that $\overline S_\theta$ is a bounded operator on $L^2(\mu, L^2(dt))$, so that in particular, $S_\theta$ is bounded in $L^2(\mu)$ and \eqref{eqTb} holds.

Finally, take $\theta(x)={\rm dist}\,\{x,\RR^d\setminus H\}$ and observe that $S_\theta$ coincides with $S$ on $\RR^d\setminus H$, so that \eqref{eqTb} implies that $S$ is bounded on $L^2(\mu |_{\CC\setminus H})$, as desired. \ep

Finally, it remains to prove Proposition~\ref{MainLemma}. This is the most technical part of the argument, which mimics the construction of Lipschitz graphs by L\'eger in \cite{Le}. We shall present the detailed argument in a full-size paper.

\bigskip

\end{document}